\documentclass{amsart}
\usepackage{amsmath,amssymb,stmaryrd}
\usepackage{amsrefs}
\usepackage{color}
\usepackage{enumerate}

\def \rr {\mathbb{R}}

\def \rn {\rr^n}
\def \sn {\mathbb{S}^{n-1}}

\def \eps {\epsilon}
\def \crits {2^\star(s)}

\newtheorem{lem}{Lemma}

\def \rn {\mathbb{R}^n}

\def \ap {\alpha_+(\gamma)}
\def \am {\alpha_-(\gamma)}

\def \dundeux {D_1^2(\rn)}
\def \eucl {\hbox{Eucl}}
\def \can {\hbox{can}}

\def \rnp {\rn\setminus\{0\}}
\def \hphi {\hat{\varphi}}
\def \hU {\hat{U}}

\newtheorem{theorem}{Theorem}[section]

\title[]{Nondegeneracy of positive solutions to nonlinear Hardy-Sobolev equations}

\author{Fr\'ed\'eric Robert}
\address{Fr\'ed\'eric Robert, Institut \'Elie Cartan, Universit\'e de Lorraine, BP 70239, F-54506 Vand{\oe}uvre-l\`es-Nancy, France}
\email{frederic.robert@univ-lorraine.fr}

\date{December 29th 2016}

\begin{document}
\begin{abstract} In this note, we prove that the kernel of the linearized equation around a positive energy solution in $\rn$, $n\geq 3$, to $-\Delta W-\gamma|x|^{-2}V=|x|^{-s}W^{\crits-1}$ is one-dimensional when $s+\gamma>0$. Here, $s\in [0,2)$, $0\leq\gamma<(n-2)^2/4$ and $\crits=2(n-s)/(n-2)$.
\end{abstract}

\maketitle

We fix $n\geq 3$, $s\in [0,2)$ and $\gamma<\frac{(n-2)^2}{4}$. We define $\crits=2(n-s)/(n-2)$. We consider a nonnegative solution $W\in C^2(\rnp)\setminus\{0\}$ to 
\begin{equation}\label{eq:V}
-\Delta W-\frac{\gamma}{|x|^{2}}W=\frac{W^{\crits-1}}{|x|^{s}}\hbox{ in }\rnp.
\end{equation}
Due to the abundance of solutions to \eqref{eq:V}, we require in addition that $W$ is an energy solution, that is $W\in \dundeux$, where $\dundeux$ is the completion of $C^\infty_c(\rn)$ for the norm $u\mapsto \Vert\nabla u\Vert_2$. Linearizing \eqref{eq:V} yields to consider
\begin{equation}\label{def:KV}
K:=\left\{\varphi\in \dundeux/\, -\Delta\varphi-\frac{\gamma}{|x|^2}\varphi=(\crits-1)\frac{W^{\crits-2}}{|x|^s}\varphi\hbox{ in }\dundeux\right\}
\end{equation}
Equation \eqref{eq:V} is conformally invariant in the following sense: for any $r>0$, define
$$W_r(x):=r^{\frac{n-2}{2}}W(rx)\hbox{ for all }x\in\rnp,$$
then, as one checks,  $W_r\in C^2(\rnp)$ is also a solution to \eqref{eq:V}, and, differentiating with respect to $r$ at $r=1$, we get that
$$-\Delta Z-\frac{\gamma}{|x|^2}Z=(\crits-1)\frac{W^{\crits-2}}{|x|^s}Z\hbox{ in }\rnp,
$$
where 
$$Z:=\frac{d}{dr}{W_r}_{|r=1}= \sum_ix^i\partial_i W+\frac{n-2}{2}W\in \dundeux.$$
Therefore, $Z\in K$. We prove that this is essentially the only element:
\begin{theorem}\label{th:main} We assume that $\gamma\geq 0$ and that $\gamma+s>0$. Then $K=\rr Z$. In other words, $K$ is one-dimensional.
\end{theorem}
Such a result is useful when performing Liapunov-Schmidt's finite dimensional reduction.
When $\gamma=s=0$, the equation \eqref{eq:V} is also invariant under the translations $x\mapsto W(x-x_0)$ for any  $x_0\in\rn$, and the kernel $K$ is of dimension $n+1$ (see Rey \cite{Rey} and also Bianchi-Egnell \cite{BE}). After this note was completed, we learnt that Dancer-Gladiali-Grossi \cite{dgg} proved Theorem \ref{th:main} in the case $s=0$, and that their proof can be extended to our case, see also Gladiali-Grossi-Neves \cite{ggn}.

\medskip\noindent This note is devoted to the proof of Theorem \ref{th:main}. Since $\gamma+s>0$, it follows from Chou-Chu \cite{ChouChu}, that there exists $r>0$ such that $W=\lambda^{\frac{1}{\crits-2}}U_r$, where
$$U(x):=\left(|x|^{\frac{2-s}{n-2}\am}+|x|^{\frac{2-s}{n-2}\ap} \right)^{-\frac{n-2}{2-s}}.$$
with
$$\eps:=\sqrt{\frac{(n-2)^2}{4}-\gamma}\hbox{  and  }\alpha_{\pm}(\gamma):=\frac{n-2}{2}\pm\sqrt{\frac{(n-2)^2}{4}-\gamma}.$$
As one checks, $U\in \dundeux\cap C^\infty(\rnp)$ and 
\begin{equation}\label{eq:U}
-\Delta U-\frac{\gamma}{|x|^2}U=\lambda\frac{U^{\crits-1}}{|x|^s}\hbox{ in }\rnp,\hbox{ with }\lambda:=4\frac{n-s}{n-2}\eps^2.
\end{equation}
Therefore, proving Theorem \ref{th:main} reduces to prove that $\tilde{K}$ is one-dimensional, where
\begin{equation}\label{def:tK}
\tilde{K}:=\left\{\varphi\in \dundeux/\, -\Delta\varphi-\frac{\gamma}{|x|^2}\varphi=(\crits-1)\lambda\frac{U^{\crits-2}}{|x|^s}\varphi\hbox{ in }\dundeux\right\}
\end{equation}

\medskip\noindent{\bf I. Conformal transformation.}\par
\noindent We let $\sn:=\{x\in\rn/\, \sum x_i^2=1\}$ be the standard $(n-1)-$dimensional sphere of $\rn$. We endow it with its canonical metric $\can$. We define 
$$\left\{\begin{array}{cccc}
\Phi: & \rr\times\sn &\mapsto &\rnp\\
&(t,\sigma) & \mapsto & e^{-t}\sigma
\end{array}\right.$$
The map $\Phi$ is a smooth conformal diffeomorphism and $\Phi^\star\eucl=e^{-2t}(dt^2+\can)$. On any Riemannian manifold $(M,g)$, we define the conformal Laplacian as $L_g:=-\Delta_g+\frac{n-2}{4(n-1)}R_g$ where $\Delta_g:=\hbox{div}_g(\nabla)$ and $R_g$ is the scalar curvature. The conformal invariance of the Laplacian reads as follows: for a metric $g'=e^{2\omega}g$ conformal to $g$ ($\omega\in C^\infty(M)$), we have that $L_{g'}u=e^{-\frac{n+2}{2}\omega}L_g(e^{\frac{n-2}{2}\omega}u)$ for all $u\in C^\infty(M)$. It follows from this invariance that for any $u\in C^\infty_c(\rnp)$, we have that
\begin{equation}\label{transfo:delta}
(-\Delta u)\circ \Phi(t,\sigma)=e^{\frac{n+2}{2}t}\left(-\partial_{tt}\hat{u}-\Delta_{\can}\hat{u}+\frac{(n-2)^2}{4}\hat{u}\right)(t,\sigma)
\end{equation}
for all $(t,\sigma)\in\rr\times \sn$, where $\hat{u}(t,\sigma):=e^{-\frac{n-2}{2}t}u(e^{-t}\sigma)$ for all $(t,\sigma)\in \rr\times\sn$. In addition, as one checks, for any $u,v\in C^\infty_c(\rnp)$, we have that
\begin{eqnarray}
\int_{\rn}(\nabla u,\nabla v)\, dx&=& \int_{\rr\times\sn}\left(\partial_t\hat{u}\partial_t\hat{v}+\left(\nabla^\prime\hat{u},\nabla^\prime\hat{v}\right)_{\can}+\frac{(n-2)^2}{4}\hat{u}\hat{v}\right)\, dt\, d\sigma\nonumber\\
&:=&B(\hat{u},\hat{v})\label{def:B}
\end{eqnarray}
where we have denoted $\nabla^\prime\hat{u}$ as the gradient on $\sn$ with respect to the $\sigma$ coordinate. We define the space $H$ as the completion of $C_c^\infty(\rr\times\sn)$ for the norm $\Vert\cdot\Vert_H:=\sqrt{B(\cdot,\cdot)}$. As one checks, $u\mapsto \hat{u}$ extends to a bijective isometry $\dundeux\to H$.

\medskip\noindent The Hardy-Sobolev inequality asserts the existence of $K(n,s,\gamma)>0$ such that $\left(\int_{\rn}\frac{|u|^{\crits}}{|x|^s}\, dx\right)^{\frac{2}{\crits}}\leq K(n,s,\gamma)\int_{\rn}\left(|\nabla u|^2-\frac{\gamma}{|x|^2}u^2\right)\, dx$ for all $u\in C^\infty_c(\rnp)$. Via the isometry $\dundeux\simeq H$, this inequality rewrites
$$\left(\int_{\rr\times \sn}|v|^{\crits}\, dt d\sigma\right)^{\frac{2}{\crits}}\leq K(n,s,\gamma)\int_{\rr\times\sn}\left((\partial_t v)^2+|\nabla^\prime v|_{\can}^2+\eps^2v^2\right)\, dtd\sigma,$$
for all $v\in H$. In particular, $v\in L^{\crits}(\rr\times\sn)$ for all $v\in H$. 

\medskip\noindent We define $H_1^2(\rr)$ (resp. $H_1^2(\sn)$) as the completion of $C^\infty_c(\rr)$ (resp. $C^\infty(\sn)$) for the norm 
$$u\mapsto \sqrt{\int_{\rr}(\dot{u}^2+u^2)\, dx}\; \left(\hbox{resp. }u\mapsto \sqrt{\int_{\sn}(|\nabla^\prime u|^2_{\can}+u^2)\, d\sigma}\right).$$
Each norm arises from a Hilbert inner product. For any $(\varphi,Y)\in C^\infty_c(\rr)\times C^\infty(\sn)$, define $\varphi\star Y\in C^\infty_c(\rr\times\sn)$ by $(\varphi\star Y)(t,\sigma):=\varphi(t)Y(\sigma)$ for all $(t,\sigma)\in\rr\times\sn$. As one checks, there exists $C>0$ such that
\begin{equation}\label{eq:star}
\Vert \varphi\star Y\Vert_H\leq C\Vert \varphi\Vert_{H_1^2(\rr)}\Vert Y\Vert_{H_1^2(\sn)}
\end{equation}
for all $(\varphi,Y)\in C^\infty_c(\rr)\times C^\infty(\sn)$. Therefore, the operator extends continuously from $H_1^2(\rr)\times H_1^2(\sn)$ to $H$, such that \eqref{eq:star} holds for all $(\varphi,Y)\in H_1^2(\rr)\times H_1^2(\sn)$.

\begin{lem}\label{lem:2} We fix $u\in C^\infty_c(\rr\times\sn)$ and $Y\in H_1^2(\sn)$. We define 
$$u_Y(t):=\int_{\sn}u(t,\sigma)Y(\sigma)\, d\sigma=\langle u(t,\cdot),Y\rangle_{L^2(\sn)}\hbox{ for all }t\in\rr.$$
Then $u_Y\in H_1^2(\rr)$. Moreover, this definition extends continuously to $u\in H$ and there exists $C>0$ such that
$$\Vert u_Y\Vert_{H_1^2(\rr)}\leq C\Vert u\Vert_H\Vert Y\Vert_{H_1^2(\sn)}\hbox{ for all }(u,Y)\in H\times H_1^2(\sn).$$
\end{lem}

\noindent{\it Proof of Lemma \ref{lem:2}:} We let $u\in C^\infty_c(\rr\times\sn)$, $Y\in H_1^2(\sn)$ and $\varphi\in C^\infty_c(\rr)$. Fubini's theorem yields:
$$\int_{\rr}\left(\partial_t u_Y\partial_t\varphi+u_Y\varphi\right)\, dt=\int_{\rr\times\sn}\left(\partial_t u\partial_t(\varphi\star Y)+u\cdot (\varphi\star Y)\right)\, dtd\sigma$$
Taking $\varphi:=u_Y$, the Cauchy-Schwartz inequality yields
\begin{eqnarray*}
&&\Vert u_Y\Vert_{H_1^2(\rr)}^2\\
&&\leq  \sqrt{\int_{\rr\times\sn}\left((\partial_t u)^2+u^2\right)dtd\sigma}
\times \sqrt{\int_{\rr\times\sn}\left((\partial_t (u_Y\star Y))^2+ (u_Y\star Y)^2\right) dtd\sigma}\\
&&\leq C\Vert u\Vert_H\Vert u_Y\star Y\Vert_H\leq C\Vert u\Vert_H\Vert u_Y\Vert_{H_1^2(\rr)}\Vert Y\Vert_{H_1^2(\sn)},
\end{eqnarray*}
and then $\Vert u_Y\Vert_{H_1^2(\rr)}\leq C\Vert u\Vert_H\Vert Y\Vert_{H_1^2(\sn)}$. The extension follows from density.\qed

\medskip\noindent{\bf II. Transformation of the problem.} We let $\varphi\in \tilde{K}$, that is 
$$-\Delta\varphi-\frac{\gamma}{|x|^2}\varphi=(\crits-1)\lambda\frac{U^{\crits-2}}{|x|^s}\varphi\hbox{ weakly in }\dundeux.$$
Since $U\in C^\infty(\rnp)$, elliptic regularity yields $\varphi\in C^\infty(\rnp)$. Moreover, the correspondance \eqref{def:B} yields
\begin{equation}\label{eq:hphi}
-\partial_{tt}\hphi-\Delta_{\can}\hphi+\eps^2\hphi=(\crits-1)\lambda \hU^{\crits-2}\hphi
\end{equation}
weakly in $H$. Note that since $\hphi,\hU\in H$ and $H$ is continuously embedded in $L^{\crits}(\rr\times\sn)$, this formulation makes sense. Since $\varphi\in C^\infty(\rnp)$, we get that $\hphi\in C^\infty(\rr\times\sn)\cap H$ and equation \eqref{eq:hphi} makes sense strongly in $\rr\times\sn$. As one checks, we have that
\begin{equation}
\hU(t,\sigma)=\left(e^{\frac{2-s}{n-2}\eps t}+e^{-\frac{2-s}{n-2}\eps t}\right)^{-\frac{n-2}{2-s}}\hbox{ for all }(t,\sigma)\in \rr\times\sn.
\end{equation}
In the sequel, we will write $\hU(t)$ for $\hU(t,\sigma)$ for $(t,\sigma)\in \rr\times\sn$.  

\medskip\noindent The eigenvalues of $-\Delta_{\can}$ on $\sn$ are 
$$0=\mu_0<n-1=\mu_1<\mu_2<....$$
We let $\mu\geq 0$ be an eigenvalue for $-\Delta_{\can}$ and we let $Y=Y_\mu\in C^\infty(\sn)$ be a corresponding eigenfunction, that is 
$$-\Delta_{\can}Y=\mu Y\hbox{ in }\sn.$$
We fix $\psi\in C^\infty_c(\rr)$ so that $\psi\star Y\in C^\infty_c(\rr\times\sn)$. Multiplying \eqref{eq:hphi} by $\psi\star Y$, integrating by parts and using Fubini's theorem yields
$$\int_{\rr}\left(\partial_{t}\hphi_Y\partial_t\psi+(\mu+\eps^2)\hphi_Y\psi\right)\, dt=\int_{\rr}(\crits-1)\lambda \hU^{\crits-2}\hphi_Y\psi\, dt,$$
where $\hphi_Y\in H_1^2(\rr)\cap C^\infty(\rr)$. Then 
\begin{equation}
A_\mu \hphi_Y=0\hbox{ with }A_\mu:=-\partial_{tt}+(\mu+\eps^2-(\crits-1)\lambda \hU^{\crits-2})
\end{equation}
where this identity holds both in the classical sense and in the weak $H_1^2(\rr)$ sense. We claim that 
\begin{equation}\label{eq:phi:0}
\hphi_Y\equiv 0\hbox{ for all eigenfunction }Y\hbox{ of }\mu\geq n-1.
\end{equation}
We prove the claim by taking inspiration from Chang-Gustafson-Nakanishi (\cite{gustaf}, Lemma 2.1). Differentiating \eqref{eq:U} with respect to $i=1,...,n$, we get that
\begin{equation}
-\Delta\partial_i U-\frac{\gamma}{|x|^2}\partial_i U-(\crits-1)\lambda\frac{U^{\crits-2}}{|x|^s}\partial_i U=-\left(\frac{2\gamma}{|x|^{4}}U+\frac{s\lambda}{|x|^{s+2}}U^{\crits-1}\right)x_i
\end{equation}
On $\rr\times\sn$, this equation reads
$$-\partial_{tt}\hat{\partial_i U}-\Delta_{\can}\hat{\partial_i U}+\left(\eps^2-(\crits-1)\lambda \hU^{\crits-2}\right)\hat{\partial_i U}=-\sigma_i e^t \left(2\gamma\hU+s\lambda \hU^{\crits-1}\right)$$
Note that $\hat{\partial_i U}=-V\star \sigma_i$, where $\sigma_i:\sn\to \rr$ is the projection on the $x_i$'s and
\begin{equation*}
V(t):=-e^{-\frac{n-2}{2}t}U^\prime(e^{-t})=e^{(1+\eps)t}\left(\ap +\am e^{2\frac{2-s}{n-2}\eps t}\right)\left(1+e^{2\frac{2-s}{n-2}\eps t}\right)^{-\frac{n-s}{2-s}}>0
\end{equation*}
for all $t\in\rr$. Since $-\Delta_{\can}\sigma_i=(n-1)\sigma_i$ (the $\sigma_i$'s form a basis of the second eigenspace of $-\Delta_{\can}$), we then get that
$$A_\mu V\geq A_{n-1}V= e^t\left(2\gamma\hU+s\lambda \hU^{\crits-1}\right)>0\hbox{ for all }\mu\geq n-1\hbox{ and }V>0.$$
Note that for $\gamma>0$, we have that $\am>0$, and that for $\gamma=0$, we have that $\am=0$. As one checks, we have that
\begin{eqnarray*}
(i)\;\left\{\left(\gamma>0\hbox{ and }\eps>1\right)\hbox{ or }\left(\gamma=0\hbox{ and }s<\frac{n}{2}\right)\right\}&\Rightarrow & V\in H_1^2(\rr)\\
(ii)\; \left\{\left(\gamma>0\hbox{ and }\eps\leq1\right)\hbox{ or }\left(\gamma=0\hbox{ and }s\geq \frac{n}{2}\right)\right\}&\Rightarrow & V\notin L^2((0,+\infty))
\end{eqnarray*}
\medskip\noindent{\it Assume that case (i) holds:} in this case, $V\in H_1^2(\rr)$ is a distributional solution to $A_\mu V>0$ in $H_1^2(\rr)$. We define $m:=\inf \{\int_{\rr}\varphi A_\mu \varphi\, dt\}$, where the infimum is taken on $\varphi\in H_1^2(\rr)$ such that $\Vert\varphi\Vert_2=1$. We claim that $m>0$. Otherwise, it follows from Lemma \ref{lem:3} below that the infimum is achieved, say by $\varphi_0\in H_1^2(\rr)\setminus \{0\}$ that is a weak solution to $A_\mu\varphi_0=m\varphi_0$ in $\rr$. Since $|\varphi_0|$ is also a minimizer, and due to the comparison principle, we can assume that $\varphi_0>0$. Using the self-adjointness of $A_\mu$, we get that $0\geq m\int_{\rr}\varphi_0V\, dt=\int_{\rr}(A_\mu \varphi_0)V\, dt=\int_{\rr}(A_\mu V)\varphi_0\, dt>0$, which is a contradiction. Then $m>0$. Since $A_\mu\varphi_Y=0$, we then get that $\varphi_Y\equiv 0$ as soon as $\mu\geq n-1$. This ends  case (i).

\medskip\noindent {\it Assume that case (ii) holds:} we assume that $\varphi_Y\not\equiv 0$. It follows from Lemma \ref{lem:4} that $V(t)=o(e^{-\alpha |t|})$ as $t\to -\infty$ for all $0<\alpha<\sqrt{\eps^2+n-1}$. As one checks with the explicit expression of $V$, this is a contradiction when $\eps<\frac{n-2}{2}$, that is when $\gamma>0$. Then we have that $\gamma=0$ and $\eps=\frac{n-2}{2}$. Since $\frac{n}{2}\leq s<2$, we have that $n=3$. As one checks, $(\mu+\eps^2-(\crits-1)\lambda \hU^{\crits-2})>0$ for $\mu\geq n-1$ as soon as $n=3$ and $s\geq 3/2$. Lemma \ref{lem:4} yields $\varphi_Y\equiv 0$, a contradiction. So $\varphi_Y\equiv 0$, this ends case (ii).

\medskip\noindent These steps above prove \eqref{eq:phi:0}. Then, for all $t\in\rr$, $\hphi(t,\cdot)$ is orthogonal to the eigenspaces of $\mu_i$, $i\geq 1$, so it is in the eigenspace of $\mu_0=0$ spanned by $1$, and therefore $\hphi=\hphi(t)$ is independent of $\sigma\in\sn$. Then 
$$-\hphi^{\prime\prime}+(\eps^2-(\crits-1)\lambda \hU^{\crits-2})\hphi=0\hbox{ in }\rr\hbox{ and }\hphi\in H_1^2(\rr).$$
It follows from Lemma \ref{lem:5} that the space of such functions is a most one-dimensional. Going back to $\varphi$, we get that $\tilde{K}$ is of dimension at most one, and then so is $K$. Since $Z\in K$, then $K$ is one dimensional and $K=\rr Z$. This proves Theorem \ref{th:main}.

\medskip\noindent{\bf III. Auxiliary lemmas.}

\begin{lem}\label{lem:5} Let $q\in C^0(\rr)$. Then 
$$\hbox{dim}_{\rr}\{\varphi\in C^2(\rr)\cap H_1^2(\rr)\hbox{ such that }-\ddot{\varphi}+q\varphi=0\}\leq 1.$$
\end{lem}
\noindent{\it Proof of Lemma \ref{lem:5}:} Let $F$ be this space. Fix $\varphi,\psi\in F\setminus\{0\}$: we prove that they are linearly dependent. Define the Wronskian $W:=\varphi \dot{\psi}-\dot{\varphi}\psi$. As one checks, $\dot{W}=0$, so $W$ is constant. Since $\varphi,\dot{\varphi},\psi,\dot{\psi}\in L^2(\rr)$, then $W\in L^1(\rr)$ and then $W\equiv 0$. Therefore, there exists $\lambda\in\rr$ such that $(\psi(0),\dot{\psi}(0))=\lambda (\varphi(0),\dot{\varphi}(0))$, and then, classical ODE theory yields $\psi=\lambda\varphi$. Then $F$ is of dimension at most one.\qed

\begin{lem}\label{lem:3} Let $q\in C^0(\rr)$ be such that there exists $A>0$ such that $\lim_{t\to\pm\infty}q(t)=A$, and define
$$m:=\inf_{\varphi\in H_1^2(\rr)\setminus\{0\}}\frac{\int_{\rr}\left(\dot{\varphi}^2+q\varphi^2\right)\, dt}{\int_{\rr}\varphi^2\, dt}.$$
Then either $m>0$, or the infimum is achieved.
\end{lem}
\noindent Note that in the case $q(t)\equiv A$, $m=A$ and the infimum is not achieved.\par
\noindent{\it Proof of Lemma \ref{lem:3}:} As one checks, $m\in\rr$ is well-defined. We let $(\varphi_i)_i\in H_1^2(\rr)$ be a minimizing sequence such that $\int_{\rr}\varphi_i^2\, dt=1$ for all $i$, that is $\int_{\rr}\left(\dot{\varphi}_i^2+q\varphi_i^2\right)\, dt=m+o(1)$ as $i\to +\infty$. Then $(\varphi_i)_i$ is bounded in $H_1^2(\rr)$, and, up to a subsequence, there exists $\varphi\in H_1^2(\rr)$ such that $\varphi_i\rightharpoonup \varphi$ weakly in $H_1^2(\rr)$ and $\varphi_i\to \varphi$ strongly in $L^2_{loc}(\rr)$ as $i\to +\infty$. We define $\theta_i:=\varphi_i-\varphi$. Since $\lim_{t\to \pm\infty}(q(t)-A)=0$ and $(\theta_i)_i$ goes to $0$ strongly in $L^2_{loc}$, we get that $\lim_{i\to +\infty}\int_{\rr}(q(t)-A)\theta_i^2\, dt=0$. Using the weak convergence to $0$ and that $(\varphi_i)_i$ is minimizing, we get that
$$\int_{\rr}\left(\dot{\varphi}^2+q\varphi^2\right)\, dt+\int_{\rr}\left(\dot{\theta}_i^2+A\theta_i^2\right)\, dt=m+o(1)\hbox{ as }i\to +\infty.$$ 
Since $1-\Vert\varphi\Vert_2^2=\Vert\theta_i\Vert_2^2+o(1)$ as $i\to +\infty$ and $\int_{\rr}\left(\dot{\varphi}^2+q\varphi^2\right)\, dt\geq m\Vert\varphi\Vert_2^2$, we get
$$m\Vert\theta_i\Vert_2^2\geq \int_{\rr}\left(\dot{\theta}_i^2+A\theta_i^2\right)\, dt+o(1)\hbox{ as }i\to +\infty.$$
If $m\leq 0$, then $\theta_i\to 0$ strongly in $H_1^2(\rr)$, and then $(\varphi_i)_i$ goes strongly to $\varphi\not\equiv 0$ in $H_1^2$, and $\varphi$ is a minimizer for $m$. This proves the lemma.\qed

\begin{lem}\label{lem:4} Let $q\in C^0(\rr)$ be such that there exists $A>0$ such that $\lim_{t\to\pm\infty}q(t)=A$ and $q$ is even. We let $\varphi\in C^2(\rr)$ be such that $-\ddot{\varphi}+q\varphi=0$ in $\rr$ and $\varphi\in H_1^2(\rr)$.
\begin{itemize}
\item If $q\geq 0$, then $\varphi\equiv 0$.
\item We assume that there exists $V\in C^2(\rr)$ such that
$$-\ddot{V}+qV>0\; ,\; V>0\hbox{ and }V\not\in L^2((0,+\infty)).$$
Then either $\varphi\equiv 0$ or $V(t)=o(e^{-\alpha |t|})$ as $t\to -\infty$ for all $0<\alpha<\sqrt{A}$.
\end{itemize}
\end{lem}
\noindent{\it Proof of Lemma \ref{lem:4}:} We assume that $\varphi\not\equiv 0$. We first assume that $q\geq 0$. By studying the monotonicity of $\varphi$ between two consecutive zeros, we get that $\varphi$ has at most one zero, and then $\ddot{\varphi}$ has constant sign around $\pm\infty$. Therefore, $\varphi$ is monoton around $\pm\infty$ and then has a limit, which is $0$ since $\varphi\in L^2(\rr)$. The contradiction follows from studying the sign of $\ddot{\varphi}$, $\varphi$. Then $\varphi\equiv 0$ and the first part of Lemma \ref{lem:4} is proved.

\medskip\noindent We now deal with the second part and we let $V\in C^2(\rr)$ be as in the statement. We define $\psi:=V^{-1}\varphi$. Then, $-\ddot{\psi}+h \dot{\psi}+Q \psi=0$ in $\rr$ with $h,Q\in C^0(\rr)$ and $Q>0$. Therefore, by studying the zeros, $\dot{\psi}$ vanishes at most once, and then $\psi(t)$ has limits as $t\to\pm\infty$.  Since $\varphi=\psi V$, $\varphi\in L^2(\rr)$ and $V\not\in L^2(0,+\infty)$, then $\lim_{t\to +\infty}\psi(t)=0$. We claim that $\lim_{t\to-\infty}\psi(t)\neq 0$. Otherwise, the limit would be $0$. Then $\psi$ would be of constant sign, say $\psi>0$. At the maximum point $t_0$ of $\psi$, the equation would yield $\ddot{\psi}(t_0)>0$, which contradicts the maximum. So the limit of $\psi$ at $-\infty$ is nonzero, and then $V(t)=O(\varphi(t))$ as $t\to-\infty$. 

\smallskip\noindent We claim that $\varphi$ is even or odd and $\varphi$ has constant sign around $+\infty$. Since $t\mapsto \varphi(-t)$ is also a solution to the ODE, it follows from Lemma \ref{lem:5} that it is a multiple of $\varphi$, and then $\varphi$ is even or odd. Since $\dot{\psi}$ changes sign at most once, then $\psi$ changes sign at most twice.  Therefore $\varphi=\psi V$ has constant sign around $+\infty$. 

\smallskip\noindent We fix $0<A'<A$ and we let $R_0>0$ such that $q(t)>A'$ for all $t\geq R_0$. Without loss of generality, we also assume that $\varphi(t)>0$ for $t\geq R_0$. We define $b(t):=C_0e^{-\sqrt{A'}t}-\varphi(t)$ for all $t\in\rr$ with $C_0:=2\varphi(R_0)e^{\sqrt{A'}R_0}$. We claim that $b(t)\geq 0$ for all $t\geq R_0$. Otherwise $\inf_{t\geq R_0}b(t)<0$, and since $\lim_{t\to +\infty}b(t)=0$ and $b(R_0)>0$, then there exists $t_1>R_0$ such that $\ddot{b}(t_1)\geq 0$ and $b(t_1)<0$. However, as one checks, the equation yields $\ddot{b}(t_1)<0$, which is a contradiction. Therefore $b(t)\geq 0$ for all $t\geq R_0$, and then $0<\varphi(t)\leq C_0e^{-\sqrt{A'}t}$ for $t\to +\infty$. Lemma \ref{lem:4} follows from this inequality, $\varphi$ even or odd, and $V(t)=O(\varphi(t))$ as $t\to-\infty$.\qed


\begin{thebibliography}{12}

\bib{BE}{article}{
   author={Bianchi, G.},
   author={Egnell, H.},
   title={A note on  Sobolev inequality},
   journal={J. Funct. Anal.},
   volume={100},
   date={1991},
   pages={18--24},
}
\bib{gustaf}{article}{
   author={Chang, S.-M.},
   author={Gustafson, S.},
   author={Nakanishi, K.},
   author={Tsai, T.-P.},
   title={Spectra of linearized operators for NLS solitary waves},
   journal={SIAM J. Math. Anal.},
   volume={39},
   date={2007/08},
   number={4},
   pages={1070--1111},
}

\bib{ChouChu}{article}{
   author={Chou, K.-S.},
   author={Chu, C.-W.},
   title={On the best constant for a weighted Sobolev-Hardy inequality},
   journal={J. London Math. Soc. (2)},
   volume={48},
   date={1993},
   number={1},
   pages={137--151},
   
}
	
	
\bib{dgg}{article}{
 author={Dancer, N.},
 author={Gladiali, F.},
 author={Grossi, M.},
 title={On the Hardy-Sobolev equation},
 journal={Proc. Roy. Soc. Edinburgh Sect. A},
 note={In press},
}

%\bib{mussowei}{article}{
%   author={Musso, M.},
%   author={Wei, J.},
%   title={Nondegeneracy of nodal solutions to the critical Yamabe problem},
%   journal={Comm. Math. Phys.},
%   volume={340},
%   date={2015},
%   number={3},
%   pages={1049--1107},
%}

%\bib{MR2975378}{article}{
%   author={Musso, M.},
%   author={Wei, J.},
%   title={Nonradial solutions to critical elliptic equations of
%   Caffarelli-Kohn-Nirenberg type},
%   journal={Int. Math. Res. Not. IMRN},
%   date={2012},
%   number={18},
%   pages={4120--4162},
% }

\bib{ggn}{article}{
   author={Gladiali, F.},
   author={Grossi, M.},
   author={Neves, S. L. N.},
   title={Nonradial solutions for the H\'enon equation in $\Bbb{R}^N$},
   journal={Adv. Math.},
   volume={249},
   date={2013},
   pages={1--36},
   
}


\bib{Rey}{article}{
   author={Rey, O.},
   title={The role of the Green's function in a nonlinear elliptic equation involving the critical Sobolev exponent},
   journal={J. Funct. Anal.},
   volume={89},
   date={1990},
   number={1},
   pages={1--52},
}

\end{thebibliography}
\end{document}